\documentclass[11pt,a4paper]{article}
\usepackage[latin1]{inputenc}
\usepackage{amsmath}
\usepackage{amsfonts}
\usepackage{amssymb}
\usepackage{makeidx}
\numberwithin{equation}{section}
\usepackage[width=12.00cm, height=19.50cm, left=2.50cm, right=2.50cm, top=2.50cm, bottom=2.cm]{geometry}
\author{Kuldeep Singh Gehlot}
\title{Two Parameter Gamma Function and It's Properties }
\begin{document}
\maketitle
\begin{center}
Government College Jodhpur,\\
JNV University Jodhpur, Rajasthan, India-306401.\\
Email: drksgehlot@rediffmail.com
\end{center}
\section*{Abstract}
In this paper we introduce the Two Parameter Gamma Function, Beta Function and Pochhammer Symbol. We named them, as p - k Gamma Function, p - k Beta Function and p - k Pochhammer Symbol and denoted as $ _{p}\Gamma_{k}(x), $  $ _{p}B_{k}(x,y) $ and $ _{p}(x)_{n,k} $ respectively. We prove the several identities for $ _{p}\Gamma_{k}(x), $  $ _{p}B_{k}(x,y) $ and $ _{p}(x)_{n,k} $ those satisfied by the classical Gamma, Beta and Pochhammer Symbol. Also we provide the integral representation for the $ _{p}\Gamma_{k}(x) $ and  $ _{p}B_{k}(x,y) $.\\\\
\textbf{Mathematics Subject Classification :} 33B15.\\\\
\textbf{Keywords:} Two Parameter Pochhammer Symbol, Two Parameter Gamma Function, Two Parameter Beta Function, Two Parameter Psi Function, p - k Hypergeometric Function. 
\section{Introduction}
The main aim of this paper is to introduce Two Parameter Pochhammer Symbol, Two Parameter Gamma Function and Two Parameter Beta Function. p - k Gamma Function is the deformation of the classical Gamma Function, such that          
$ _{p}\Gamma_{k}(x) \Rightarrow\: _{k}\Gamma_{k}(x) = \Gamma_{k}(x)  $ as $ p = k $ and $ _{p}\Gamma_{k}(x) \Rightarrow \:_{1}\Gamma_{1}(x)= \Gamma (x)  $ as $ p,k\rightarrow 1 $.\\
In section 2, we defined Two Parameter Pochhammer Symbol denoted as $ _{p}(x)_{n,k} $, with its convergent conditions. Two Parameter Pochhammer Symbol is the deformation of the classical Pochhammer Symbol, such that $ _{p}(x)_{n,k} \Rightarrow  \: _{k}(x)_{n,k} = (x)_{n,k}  $ as $ p=k $ and  $ _{p}(x)_{n,k} \Rightarrow  \:  _{1}(x)_{n,1} = (x)_{n}  $ as $ p=k=1. $  Also we derived two parameter Pochhammer symbol  in terms of the elementary symmetric function, and evaluate it's derivative identities. It is most natural to relate the two parameter Pochhammer symbol to the two parameter Gamma Function is defined. We evaluate integral representation of two parameter Gamma Function, also represent two parameter Gamma Function  into different infinite product forms and so many recurrence relations are evaluated. \\
In section 3, we defined two parameter Beta Function and two parameter Psi Function. Also evaluate some recurrence relations and functional relation with classical Beta Function.\\
Section 4, deal with the definition of Hypergeometric function with Two Parameter Pochhammer Symbol, known as p-k Hypergeometric function. Also we evaluate the Differential Equation, Functional relation with Classical Hypergeometric function and Integral Representation of p-k Hypergeometric function.\\
Throughout this paper Let $ C,R^{+}, Re(),Z^{-} and N $ be the sets of complex numbers, positive real numbers, real part of complex number, negative integer and natural numbers respectively. 
\section{p - k Pochhammer Symbol and p - k Gamma Function }
In this section we introduce p - k Pochhammer Symbol and p - k Gamma Function. We evaluate $ _{p}\Gamma_{k}(x) $ in terms of limit, recurrence formulas and infinite products. 
\subsection{Definition}\label{1.1}
Let $ x\in C ;  k,p \in R^{+}-\lbrace 0 \rbrace $ and $ Re(x)>0, n\in N, $ the p - k Pochhammer Symbol (i.e. Two Parameter Pochhammer Symbol), $ _{p}(x)_{n,k} $ is given by 
\begin{equation}
_{p}(x)_{n,k}=(\frac{xp}{k})(\frac{xp}{k}+p)(\frac{xp}{k}+2p).........(\frac{xp}{k}+(n-1)p).
\end{equation}
For $ s,n \in N  $ with $ 0\leq s\leq n, $ the $ s^{th} $ elementary symmetric function 
\[e^{n}_{s}(x_{1},x_{2},.....,x_{n})= \sum_{{1\leq i_{1} \leq i_{2} \leq... \leq i_{s} \leq n}} x_{i_{1}}.....x_{i_{s}}\]
on the variables $ x_{1},x_{2},.....,x_{n}.$\\\\ 
\textbf{Theorem 2.1} Formula for the p - k Pochhammer Symbol (i.e. Two Parameter Pochhammer Symbol) in terms of the elementary symmetric function is given by 
\begin{equation}
_{p}(x)_{n,k}=\sum_{s=0}^{n-1}p^{n}e_{s}^{n-1}(1,2,...,(n-1))(\frac{x}{k})^{n-s}.
\end{equation}
Where $ x\in C ;  k,p \in R^{+}-\lbrace 0 \rbrace $ and $ Re(x)>0, n\in N. $\\
Proof: The well known identity for elementary symmetric polynomials appear when expand a linear factorization of a monic polynomial
\begin{equation}
\prod_{j=1}^{n}(\lambda + X_{j})= \lambda^{n}e_{0}^{n-1}+\lambda^{n-1}e_{1}^{n-1}+.....+\lambda e_{n-1}^{n-1} =\sum_{s=0}^{n-1}e_{s}^{n-1}(1,2,...,(n-1))(\lambda)^{n-s}.
\end{equation} 
Using equation (2.1), we have the desired result.\\\\
\textbf{Theorem 2.2} The derivative identities for p - k Pochhammer Symbol. 
\begin{equation}
\frac{\partial}{\partial k}[_{p}(x)_{n,k}]= -\frac{n}{k}\:_{p}(x)_{n,k} + \frac{p}{k}\sum_{s=1}^{n-1}s\:_{p}(x)_{n,k}\: _{p}(x+(s+1)k)_{n-1-s,k}.
\end{equation}
\begin{equation}
\frac{\partial}{\partial p}[_{p}(x)_{n,k}]= \frac{n}{p}\:_{p}(x)_{n,k}.
\end{equation}
Where $ x\in C ;  k,p \in R^{+}-\lbrace 0 \rbrace $ and $ Re(x)>0, n\in N. $\\
Proof: Using the definition (2.1) and logarithmic derivatives, we have the desired result.
\subsection{Definition}
For $ x\in C/kZ^{-};  k,p \in R^{+}-\lbrace 0 \rbrace $ and $ Re(x)>0, n\in N, $ the p - k Gamma Function (i.e. Two Parameter Gamma Function), $_{p}\Gamma_{k}(x)$ is given by
\begin{equation}
   _{p}\Gamma_{k}(x)=\frac{1}{k}\lim_{n\rightarrow \infty} \dfrac{n!p^{n+1}(np)^{\frac{x}{k}}}{_{p}(x)_{n+1,k}}.
\end{equation}  
 or
\begin{equation}
   _{p}\Gamma_{k}(x)=\frac{1}{k}\lim_{n\rightarrow \infty} \dfrac{n!p^{n+1}(np)^{\frac{x}{k}-1}}{_{p}(x)_{n,k}}.
\end{equation} 
\textbf{Theorem 2.3} Given $ x\in C/kZ^{-} ;  k,p,s,r \in R^{+}-\lbrace 0 \rbrace $ and $ Re(x)>0, n\in N, $ the following identities holds,
\begin{equation}
 _{p}(x)_{n,s}= \:_{p}(\frac{kx}{s})_{n,k}.
\end{equation}
\begin{equation}
 _{p}(x)_{n,s}= (\frac{p}{s})^{n} \; _{s}(\frac{kx}{s})_{n,k}.
\end{equation}
\begin{equation}
 _{p}(x)_{n,k}= (\frac{p}{s})^{n}\:_{s}(x)_{n,k}.
\end{equation}
\begin{equation}
 _{p}\Gamma_{s}(x)= \frac{k}{s} \; _{p}\Gamma_{k}(\frac{kx}{s}).
 \end{equation} 
\begin{equation}
 _{r}\Gamma_{s}(x)= \frac{k}{s} \;(\frac{r}{p})^{\frac{x}{s}}\; _{p}\Gamma_{k}(\frac{kx}{s}).
 \end{equation} 
 \begin{equation}
 _{r}\Gamma_{k}(x)= (\frac{r}{p})^{\frac{x}{k}}\; _{p}\Gamma_{k}(x).
 \end{equation}
 Proof: Property (2.8), (2.9), (2.10) follows directly from definition (2.1) and the results (2.11), (2.12), (2.13) will follow directly  by using equation (2.6).\\\\
\textbf{Theorem 2.4} Given $ x\in C / kZ^{-}; k,p\in R^{+}-\lbrace 0 \rbrace $ and $ Re(x)>0, $
then the integral representation of p - k Gamma Function is given by
\begin{equation}
_{p}\Gamma_{k}(x)=\int^{\infty}_{0}e^{-\frac{t^{k}}{p}}t^{x-1}dt.
 \end{equation} 
Proof: Consider the right hand side integral and ([3],  Page 2) Tannery's Theorem and equation (2.7), we have\\
\[\int^{\infty}_{0}e^{-\frac{t^{k}}{p}}t^{x-1}dt=\lim _{n \rightarrow \infty}\int^{(np)^{\frac{1}{k}}}_{0}(1-\frac{t^{k}}{np})^{n}\:t^{x-1}dt.\]
Let $ A_{n,i}(x), i = 0,1,...,n, $ be given by
\[A_{n,i}(x) = \int^{(np)^{\frac{1}{k}}}_{0}(1-\frac{t^{k}}{np})^{i}\:t^{x-1}dt.\] 
Integration by parts we have the following recurrence formula,
\[A_{n,i}(x)= \frac{ki}{pxn} A_{n,i-1}(x+k).\]
Also,
\[A_{n,0}(x)= \int^{(np)^{\frac{1}{k}}}_{0}t^{x-1}dt=\frac{(np)^{\frac{x}{k}}}{x}.\]
Therefor,
\[A_{n,n}(x)= \frac{1}{k} \dfrac{n!p^{n+1}(np)^{\frac{x}{k}-1}}{_{p}(x)_{n,k}(1+\frac{x}{kn})}.\]
\[ _{p}\Gamma_{k}(x)= \lim_{n\rightarrow \infty} A_{n,n}(x) =\frac{1}{k}\lim_{n\rightarrow \infty} \dfrac{n!p^{n+1}(np)^{\frac{x}{k}-1}}{_{p}(x)_{n,k}}.\]
Which complete the proof.\\\\ 
\textbf{Theorem 2.5} Given $ x\in C / kZ^{-}; k,p\in R^{+}-\lbrace 0 \rbrace $ and $ Re(x)>0, $
then we have,
\begin{equation}
_{p}\Gamma_{k}(x)=\frac{p^{\frac{x}{k}}}{k}\prod_{n=1}^{\infty}[(1+\frac{1}{n})^{\frac{x}{k}}(1+\frac{x}{nk})^{-1}].
\end{equation}
Proof: Using equation (2.1) and (2.7), we immediately get the desire result.\\\\
\textbf{Theorem 2.6} Given $ x\in C / kZ^{-}; k,p\in R^{+}-\lbrace 0 \rbrace $ and $ Re(x)>0, $
then we have,
\begin{equation}
\frac{1}{_{p}\Gamma_{k}(x)}=\frac{x}{kp^{\frac{x}{k}}}\lim _{n\rightarrow\infty}[n^{-\frac{x}{k}}\prod_{r=1}^{n}(1+\frac{x}{rk})].
\end{equation}
Proof: Using equation (2.1) and (2.7), we immediately get the desire result.\\\\
\textbf{Theorem 2.7} Given $ x\in C / kZ^{-}; k,p\in R^{+}-\lbrace 0 \rbrace $ and $ Re(x)>0, $
then we have,
\begin{equation}
_{p}\Gamma_{k}(x)=a^{\frac{x}{k}}\int^{\infty}_{0}e^{-\frac{t^{k}}{p}a}t^{x-1}dt.
\end{equation}
Proof: Using equation (2.14), we immediately get the desire result.\\\\
\textbf{Theorem 2.8} Given $ x\in C / kZ^{-}; k,p\in R^{+}-\lbrace 0 \rbrace $ and $ Re(x)>0, $
then we have,
\begin{equation}
\frac{1}{_{p}\Gamma_{k}(x)}=\frac{x}{kp^{\frac{x}{k}}}e^{\frac{x}{k}\gamma}\prod_{n=1}^{\infty}[(1+\frac{x}{nk})e^{-\frac{x}{nk}}].
\end{equation}
Where\[ \gamma = \lim_{n \rightarrow \infty}[1+\frac{1}{2}+....+\frac{1}{n}-\log n],\] is Euler's constant.\\
Proof: Using equation (2.16), we immediately get the desire result.\\\\
\textbf{Theorem 2.9} The relation between p - k Gamma Function, k-Gamma Function and classical Gamma Function is given by,
\begin{equation}
 _{p}\Gamma_{k}(x)=(\frac{p}{k})^{\frac{x}{k}}\Gamma_{k}(x)= \frac{p^{\frac{x}{k}}}{k}\Gamma(\frac{x}{k}).
\end{equation}
Where $ x\in C / kZ^{-}; k,p\in R^{+}-\lbrace 0 \rbrace $ and $ Re(x)>0, $\\
Proof: Using (2.14) and Proposition 4, page 3 of [1], we get the desire result.\\\\
\textbf{Theorem 2.10} For $ x\in C / kZ^{-}; n,q \in N ;k,p\in R^{+}-\lbrace 0 \rbrace $ and $ Re(x)>0, $ then the relation between p - k Pochhammer Symbol, k-Pochhammer Symbol and classical Pochhammer Symbol is given by,
\begin{equation}
_{p}(x)_{n,k}=(\frac{p}{k})^n(x)_{n,k} = (p)^n(\frac{x}{k})_{n}.
\end{equation}
Also for Generalized p - k Pochhammer Symbol, we have  
\begin{equation}
_{p}(x)_{nq ,k}=(\frac{p}{k})^{nq} \:(x)_{nq,k} = (p)^{nq} \:(\frac{x}{k})_{nq}=(pq)^{nq}\:\prod_{r=1}^{q}(\dfrac{\frac{x}{k}+r-1}{q})_{n}.
\end{equation}
 \\
Proof: Using (2.1), (2.19) and Proposition 4, page 3 of [1], we get the desire result.\\\\
\textbf{Theorem 2.11} For $ x\in C/kZ^{-} ;  k,p \in R^{+}-\lbrace 0 \rbrace $ and $ Re(x)>0, n\in N. $ The fundamental equations satisfied by p - k Gamma Function, $ _{p}\Gamma_{k}(x) $ are,
\begin{equation}
_{p}(x)_{n,k}=\frac{ _{p}\Gamma_{k}(x+nk)}{ _{p}\Gamma_{k}(x)}.
\end{equation}
\begin{equation}
 _{p}\Gamma_{k}(x+k)=\frac{xp}{k}\: _{p}\Gamma_{k}(x).
\end{equation}
\begin{equation}
 _{p}\Gamma_{k}(x+nk)= p^{n}(\frac{x}{k})(\frac{x}{k}+1)......(\frac{x}{k}+(n-1))\:_{p}\Gamma_{k}(x).
\end{equation}
\begin{equation}
\frac{_{p}\Gamma_{k}(x)}{_{p}\Gamma_{k}(x-nk)}=\frac{p^{n}}{k^{n}}(x-k)(x-2k).....(x-nk).
\end{equation}
\begin{equation}
\frac{_{p}\Gamma_{k}(x)}{_{p}\Gamma_{k}(x-nk)}=(-1)^{n}\:\frac{_{p}\Gamma_{k}(-x+nk+k)}{_{p}\Gamma_{k}(-x+k)}.
\end{equation}
\begin{equation}
_{p}\Gamma_{k}(1)=\frac{p^{\frac{1}{k}}}{k}\Gamma(\frac{1}{k}).
\end{equation}
\begin{equation}
_{p}\Gamma_{k}(k)=\frac{p}{k}.
\end{equation}
\begin{equation}
_{p}\Gamma_{k}(p)=\frac{p^{\frac{p}{k}}}{k}\Gamma(\frac{p}{k}).
\end{equation}
\begin{equation}
_{p}\Gamma_{k}(x)\:_{p}\Gamma_{k}(-x)=\frac{\pi}{xk}\:\frac{1}{\sin (\frac{\pi x}{k})}.
\end{equation}
\begin{equation}
_{p}\Gamma_{k}(x)\:_{p}\Gamma_{k}(k-x)=\frac{p}{k^{2}}\:\frac{\pi}{\sin (\frac{\pi x}{k}) }.
\end{equation}
\begin{equation}
\:\prod_{r=0}^{m-1}\:_{p}\Gamma_{k}(x+\frac{kr}{m})=\frac{p^{\frac{m-1}{2}}}{k^{m-1}}(2 \pi )^{\frac{(m-1)}{2}}m^{\frac{1}{2}-\frac{mx}{k}}\:_{p}\Gamma_{k}(mx); m=2,3,4,... \: .
\end{equation}
Proof: All the results follow directly from using equation (2.1), (2.7) and (2.14). \\\\
\textbf{Theorem 2.12} For $ x\in C/kZ^{-} ;  k,p \in R^{+}-\lbrace 0 \rbrace $ and $ Re(x)>0, n\in N. $
Then the recurrence relations for p - k Pochhammer Symbol are given by,
\begin{equation}
n\:_{p}(x)_{n-1,k}=\:_{p}(x)_{n,k}- \:_{p}(x-k)_{n,k}.
\end{equation}
And
\begin{equation}
_{p}(x)_{n+j,k}=\:_{p}(x)_{j,k}\times\:_{p}(x+jk)_{n,k}.
\end{equation}
Proof: Using equation (2.20) and basic relations $ n(x)_{n-1}=(x)_{n}-(x-1)_{n}, (x)_{n+j}=(x)_{j}(x+j)_{n},$ we get the desired result.  
\section{p - k Beta Function and p - k Psi Function}
In this section, we introduce the p - k Beta Function $ _{p}B_{k}(x,y) $ and p - k Psi Function $ _{p}\psi_{k}(x,y). $ We evaluate explicit formula that relate the $ _{p}B_{k}(x,y) $ and $ _{p}\psi_{k}(x) $ to classical Beta function $ B(x,y) $ and Classical Psi function respectively $ \psi(x). $ Also prove some identities.
\subsection{Definition}
The p - k Beta Function $ _{p}B_{k}(x,y) $ is given by 
\begin{equation}
    _{p}B_{k}(x,y)=\frac{_{p}\Gamma _{k}(x)\: _{p}\Gamma _{k}(y)}{_{p}\Gamma _{k}(x+y)}; Re(x)>0, Re(y)>0.
\end{equation}  
\textbf{Theorem 3.1} The $ _{p}B_{k}(x,y) $ function satisfies the following identities. 
\begin{equation}
  _{p}B_{k}(x,y)=\frac{1}{k} \int_{0}^{1}t^{\frac{x}{k}-1}(1-t)^{\frac{y}{k}-1}dt.
 \end{equation} 
\begin{equation}
  _{p}B_{k}(x,y)=\frac{1}{k} \int_{0}^{1}\frac{t^{\frac{x}{k}-1}+t^{\frac{y}{k}-1}}{(t+1)^{\frac{x+y}{k}}}dt.
 \end{equation}
\begin{equation}
  _{p}B_{k}(x,y)= \int_{0}^{\infty}t^{x-1}(1+t^{k})^{-\frac{x+y}{k}}dt.
 \end{equation} 
\begin{equation}
  _{p}B_{k}(x,y)=\frac{1}{k}\: B(\frac{x}{k},\frac{y}{k}).
 \end{equation}
 Proof: Using the definition (3.1), we have immediately above results.
\subsection{Definition}
The logarithmic derivative of the p - k Gamma Function is known as p - k Psi Function, $ _{p}\psi_{k}(x).$
\begin{equation}
   _{p}\psi_{k}(x)=\frac{d}{dx} \ln [_{p}\Gamma _{k}(x)]=\frac{1}{_{p}\Gamma _{k}(x)}\frac{d}{dx}[_{p}\Gamma _{k}(x)].
\end{equation} 
\begin{equation}
\ln [_{p}\Gamma _{k}(x)]=\int_{1}^{x}\:_{p}\psi_{k}(x)dx.
\end{equation}
\textbf{Theorem 3.2} Some properties of $ _{p}\psi_{k}(x) $ are given by
\begin{equation}
  _{p}\psi_{k}(x)=\frac{\ln p}{k}+\psi(\frac{x}{k}).
  \end{equation}  
\begin{equation}
  _{p}\psi_{k}(x)=\frac{\ln p}{k}-\gamma- \frac{k}{x}+ x\sum_{n=1}^{\infty}\frac{1}{n(x+nk)}.
  \end{equation}   
  \begin{equation}
  _{p}\psi_{k}(x)=\frac{\ln p}{k}-\gamma + (x-k)\sum_{n=0}^{\infty}\frac{1}{(n+1)(x+nk)}.
  \end{equation} 
  Where $ \gamma  $ is Euler's Constant and $ \psi(x) $ is Classical Psi Function.\\\\
 Proof: Using the definition (3.2), we have immediately above results.\\\\
 \textbf{Theorem 3.3} The $ r^{th} $ derivative of p - k Psi Function, $ _{p}\psi_{k}(x)$ yields the result in terms of k-Zeta Function, $ \zeta_{k}(x,r), $
 \begin{equation}
 \frac{d^{r}}{dx^{r}}[\ln [_{p}\Gamma _{k}(x)]]=\frac{d^{r-1}}{dx^{r-1}}\:_{p}\psi_{k}(x)=(-1)^{r}k(r-1)!\sum_{n=0}^{\infty}\frac{1}{(x+nk)^{r}}, for \: \:r\geq 2
 \end{equation}
 Where k-Zeta Function given by definition 15, page 8 of [1]. 
 \[ \zeta_{k}(x,r)=\sum_{n=0}^{\infty}\frac{1}{(x+nk)^{r}}.\]
 Proof: Using the definition (3.2) and differentiate, we get the desired result.
\section{Hypergeometric Function}
In this section we define the Hypergeometric Function using p - k Pochhammer Symbols. Here we are use the notation of  [2].
\subsection{Definition}
Given $ x \in C, $ $ a\in C^{r}; k,p \in (R^{+})^{r}; s,t \in (R^{+})^{q}, b=(b_{1},b_{2},...,b_{q})\in C^{q} $ such that $b_{i} \in C/ s_{i}Z^{-}.$ The p-k hypergeometric function $ F(a,p,k;b,t,s;x)$ is given by 
 \begin{equation}
   F(a,p,k;b,t,s;x)= \sum_{n=0}^{\infty}\frac{\prod_{i=1}^{r}\:_{p_{i}}(a_{i})_{n,k_{i}}}{\prod_{j=1}^{q}\:_{t_{j}}(b_{j})_{n,s_{j}}}\frac{x^{n}}{n!}.
   \end{equation} 
By using Ratio Test we can show that the series (4.1) converges for all finite $ x $ if $ r\leq q.$ If $ r > q+1, $ the series diverges and if $ r = q+1, $ it converges for all $ x $ such that $ \vert x \vert < \vert\frac{t_{1}t_{2}.....t_{q}}{p_{1}p_{2}.....p_{r}}\vert. $ \\\\
\textbf{Theorem 4.1} Given $ x \in C, $ $ a\in C^{r}; k,p \in (R^{+})^{r}; s,t \in (R^{+})^{q}, b=(b_{1},b_{2},...,b_{q})\in C^{q} $ such that $b_{i} \in C/ s_{i}Z^{-}.$ Then the Functional relation between p - k Hypergeometric Function and Classical Hypergeometric Function is given by,
\begin{equation}
F(a,p,k;b,t,s;x)=F(\frac{a}{k};\frac{b}{s};\frac{\prod_{i=1}^{r}p_{i}}{\prod_{j=1}^{q}t_{j}}x).
\end{equation} 
Proof: Using definition (2.20), we get above result.\\\\
\textbf{Theorem 4.2} The Differential Equation of p - k Hypergeometric Function is given by 
\begin{equation}
[\theta \prod_{j=1}^{q}(\theta + \frac{b_{j}}{s_{j}}-1)-Ax\prod_{i=1}^{r}(\theta + \frac{a_{i}}{k_{i}}) ]W=0.
\end{equation}
Where $ \theta = x\frac{d}{dx}, $ $ A=\dfrac{\prod_{i=1}^{r}p_{i}}{\prod_{j=1}^{q}t_{j}} $ and $ W= F(a,p,k;b,t,s;x).$\\
For $ r\leq q+1,$ $i=1,2,...,r $ and $ j=1,2,...,q $ when no $ \frac{b_{j}}{s_{j}} $ is a negative integer or zero and no two $ \frac{b_{j}}{s_{j}} $ is differ by an integer or zero.\\
Proof: Using Function relation (4.2), we get the desired result.\\\\
\textbf{Theorem 4.3} For any $ a\in C; k,p > 0 $ and $ \vert x \vert <\frac{1}{p}, $ the following identity holds
\begin{equation}
\sum_{n=0}^{\infty}\frac{_{p}(a)_{n,k}\: x^{n}}{n!}= (1-xp)^{-\frac{a}{k}}. 
\end{equation}
Proof: Using (2.20), we get immediately the desired result.\\\\
\textbf{Theorem 4.4} Given $ x \in C, $ $ a\in C^{r}; k,p \in (R^{+})^{r}; s,t \in (R^{+})^{q}, b=(b_{1},b_{2},...,b_{q})\in C^{q} $ such that $b_{i} \in C/ s_{i}Z^{-}.$ The Integral Representation of 
p - k Hypergeometric Function is given by,
\begin{equation}
F(a,p,k;b,t,s;x)=\prod_{i=1}^{r}\prod_{j=1}^{q}\frac{\Gamma(\frac{b_{j}}{s_{j}})}{\Gamma(\frac{a_{i}}{k_{i}})\Gamma(\frac{b_{j}}{s_{j}}-\frac{a_{i}}{k_{i}})}\int_{0}^{1}t^{\frac{a_{i}}{k_{i}}-1}(1-t)^{\frac{b_{j}}{s_{j}}-\frac{a_{i}}{k_{i}}-1}e^{\frac{p_{i}}{t_{j}}xt}dt.
 \end{equation} 
 Proof:Using (4.2), we get immediately the desired result.\\
\section*{References}
\label{1}[1]	Diaz, R. and Pariguan, E. On hypergeometric functions and Pochhammer k-symbol. Divulgaciones Mathematicas, Vol. 15 No. 2 (2007) 179-192.\\
\label{2}[2]  Earl D. Rainville, Special Function, The Macmillan Company, New york,1963.\\
\label{3}[3]  Erdelyi, A., Higher Transcendental  Function Vol. 1, McGraw-Hill Book Company, New York, 1953. 
\end{document}